\documentclass[11pt]{article}
\RequirePackage[OT1]{fontenc}
\RequirePackage[amsthm,amsmath,natbib]{imsart}
\usepackage{amsfonts,graphicx,url}
\newcommand {\apgt} {\ {\raise-.5ex\hbox{$\buildrel>\over\sim$}}\ }
\newcommand {\aplt} {\ {\raise-.5ex\hbox{$\buildrel<\over\sim$}}\ }
\usepackage{multirow}
\usepackage{hhline}

\begin{document}


\begin{frontmatter}

\title{Moment Consistency of the Exchangeably Weighted Bootstrap for
Semiparametric M-Estimation}
\author{Guang Cheng\thanksref{t1}}
\runauthor{Guang Cheng} \runtitle{Semiparametric Bootstrap
Moment Consistency}

Purdue University

\thankstext{t1}{Associate Professor, Department of Statistics, Purdue University,
West Lafayette, IN 47907-2066, Email: chengg@purdue.edu. Research
sponsored by NSF (DMS-0906497, CAREER Award DMS-1151692, DMS-1418042), Simons
Foundation 305266. Guang Cheng was on sabbatical at Princeton while the finalization of this work was carried out; he would like to thank the Princeton ORFE department for its hospitality and support.}

\maketitle
\begin{abstract}
The bootstrap variance estimate is widely used in semiparametric
inferences. However, its theoretical validity is a well known open
problem. In this paper, we provide a {\em first} theoretical study
on the bootstrap moment estimates in semiparametric models.
Specifically, we establish the bootstrap moment consistency of the
Euclidean parameter which immediately implies the consistency of
$t$-type bootstrap confidence set. It is worth pointing out that
the only additional cost to achieve the bootstrap moment
consistency in contrast with the distribution consistency is to simply strengthen the
$L_1$ maximal inequality condition required in the latter to the
$L_p$ maximal inequality condition for $p\geq 1$. The general
$L_p$ multiplier inequality developed in this paper is also of independent interest. These general
conclusions hold for the bootstrap methods with exchangeable
bootstrap weights, e.g., nonparametric bootstrap and Bayesian bootstrap. Our general theory is illustrated in the
celebrated Cox regression model.
\end{abstract}


\end{frontmatter}

\newtheorem{theorem}{\indent \sc Theorem}
\newtheorem{corollary}{\indent \sc Corollary}
\newtheorem{lemma}{\indent \sc Lemma}
\newtheorem{proposition}{\indent \sc Proposition}
\newtheorem{remark}{\indent \sc Remark}
\newcommand{\phif}{\textsc{igf}}
\newcommand{\sign}{\mbox{sign}}
\newcommand{\phgf}{\textsc{gf}}
\newcommand{\fix}{$\textsc{gf}_0$}
\newcommand{\mb}[1]{\mbox{\bf #1}}
\newcommand{\Exp}[1]{\mbox{E}\left[#1\right]}
\newcommand{\pr}[1]{\mbox{P}\left[#1\right]}
\newcommand{\pp}[0]{\mathbb{P}}
\newcommand{\ee}[0]{\mbox{E}}
\newcommand{\re}[0]{\mathbb{R}}
\newcommand{\argmax}[0]{\mbox{argmax}}
\newcommand{\argmin}[0]{\mbox{argmin}}
\newcommand{\ind}[0]{\mbox{\Large\bf 1}}
\newcommand{\narrow}{\stackrel{n\rightarrow\infty}{\longrightarrow}
}
\newcommand{\weakpn}{\stackrel{P_n}{\leadsto}}
\newcommand{\weakpnboot}{\mbox{\raisebox{-1.5ex}{$\stackrel
{\mbox{\scriptsize $P_n$}}{\stackrel{\mbox{\normalsize
$\leadsto$}} {\mbox{\normalsize $\circ$}}}$}}\,}
\newcommand{\ol}[1]{\overline{#1}}
\newcommand{\avgse}[1] { \bar{\hat{\Sigma}}_{#1} }
\newcommand{\mcse}[1]  { \Sigma^{*}_{#1} }
\newcommand{\po}{\textsc{po}}
\newcommand{\ph}{\textsc{ph}}
\newcommand{\x}{\mathbf{x}}
\newcommand{\xw}{\mathbf{xw}}

\def\boxit#1{\vbox{\hrule\hbox{\vrule\kern6pt
          \vbox{\kern6pt#1\kern6pt}\kern6pt\vrule}\hrule}}
\def\jhcomment#1{\vskip 2mm\boxit{\vskip 2mm{\color{red}\bf#1}
 {\color{blue}\bf -- JH\vskip 2mm}}\vskip 2mm}

\def\jhcommentin#1{{(\color{red}\bf{#1}{\color{blue}\bf -- JH})}}



\def\boxit#1{\vbox{\hrule\hbox{\vrule\kern6pt
          \vbox{\kern6pt#1\kern6pt}\kern6pt\vrule}\hrule}}

\def\chengcomment#1{\vskip 2mm\boxit{\vskip 2mm{\color{blue}\bf#1}
{\color{red}\bf -- Cheng\vskip 2mm}}\vskip 2mm}

\def\chengcommentin#1{{(\color{blue}\bf{#1}{\color{red}\bf -- Cheng}})}


\section{Introduction}
In semiparametric models, the asymptotic variance estimate for the
Euclidean parameter is required in the construction of confidence
sets and test statistics based on the asymptotic normality result.
For example, in the bootstrap inferences, the asymptotic variance
estimate is needed to build the t-type confidence set which is
known to have smaller coverage probability error than the
percentile/hybrid confidence sets; see \citep{st96}. In general,
the explicit variance estimation is not feasible due to the
presence of an infinite dimensional nuisance parameter; see
\citep{bkrw98, vw96} for numerous examples. In the literature,
there are two existing estimation procedures, i.e., the profile
sampler \citep{lkf05} and the observed profile information
\citep{mv99}. The former (latter) method requires a careful choice
of the prior on the Euclidean parameter (of the step size in
calculating discretized information estimate). Subsampling
\citep{pr94b} is another possibility, but the optimal subsample
size is difficult to choose in practice. In contrast, the
bootstrap can estimate the asymptotic variance without involving
any tuning parameter, and thus becomes one standard
semiparametric inference procedure. Various types of bootstrap
variance estimate based on different sampling schemes such as
nonparametric bootstrap or weighted bootstrap are proposed in a
broad class of semiparametric models ranging from the simple
partly linear models (Chapter 2 of \citep{hlg00}), to the
complicated proportional hazards frailty regression models
(\citep{klf04}), and widely used semiparametric conditional moment
models (\citep{cp09}). More examples can be found in Kosorok
(2008). However, the theoretical validity of the bootstrap
variance estimate is a well known open problem.

Cheng and Huang (2010) have recently proven that the exchangeably
weighted bootstrap is asymptotically consistent in estimating the
distribution of the M-estimate of Euclidean parameter. However,
this distributional consistency does not imply the consistency of
the bootstrap variance estimators. Nishiyama (2010) and Kato
(2011) have shown the moment convergence of the (nonparametric
bootstrap) M-estimate in parametric models. Inspired by these
recent developments, we provide a first theoretical study on the
bootstrap moment estimates in semiparametric models. Specifically,
we establish the bootstrap moment consistency of the Euclidean
parameter which immediately implies the consistency of $t$-type
bootstrap confidence set with the help of the conditional
Slutsky's Lemma. It is worthy pointing out that the only
additional cost to achieve the bootstrap moment consistency in
contrast with the distribution consistency is to simply strengthen
the $L_1$ maximal inequality condition required in the latter to
the $L_p$ maximal inequality condition for $p\geq 1$. The general
$L_p$ multiplier inequality developed in this paper is the key
technical tool, and is also of independent interest. Our general
conclusions hold for the bootstrap methods with exchangeable
bootstrap weights, e.g., nonparametric bootstrap, and apply to a
broad class of semiparametric models with root-n convergent
nuisance parameters, e.g., Cox regression model, proportional odds
model and case control studies with a missing covariate
\citep{mv01}. The classical Cox regression model is used to
illustrate the practicality of the required conditions. Some simulations studies are also conducted for this model. As far as
we are aware, this paper presents the first theoretical studies on
the bootstrap variance consistency in semiparametric models.

\section{Preliminary}
\subsection{Semiparametric M-Estimation}
The semiparametric M-estimation, including the maximum likelihood
estimation as a special case, refers to a general method of
estimation. Let $\theta\in\Theta\subset\mathbb {R}^d$ be a
Euclidean parameter of interest and $\eta\in\mathcal{H}$ be an
infinite dimensional nuisance parameter with the norm $d(\cdot)$.
The semiparametric M-estimator $(\widehat{\theta},\widehat{\eta})$
is obtained by optimizing some objective function $m(\theta,\eta)$
based on the observations $(X_1,\ldots,X_n)$:
\begin{eqnarray}
{\textstyle (\widehat{\theta},\widehat{\eta})
=\arg\sup_{\theta\in\Theta,\eta\in\mathcal{H}}\sum_{i=1}^{n}
m(\theta,\eta)(X_{i}).}\label{mestimate}\vspace{-0.15in}
\end{eqnarray}
\noindent The form of the objective function depends on the
context. For example, it could be the log-likelihood,
quasi-likelihood \citep{mvg97} or some pseudo-likelihood function,
e.g., \citep{wz07}. Define
$(\theta_0,\eta_0)=\arg\sup_{\theta\in\Theta,\eta\in\mathcal{H}}E_Xm(\theta,\eta)(X)$.
Under mild conditions, Cheng and Huang (2010) show that
\begin{eqnarray}
\sqrt{n}(\widehat\theta-\theta_0)\overset{d}{\longrightarrow} N(0,
\Sigma).\label{mestdis}
\end{eqnarray}
Note that $\widehat\theta$ is semiparametric efficient and
$\Sigma$ is the inverse of the efficient information matrix when
$m(\theta,\eta)$ is the log-likelihood function.

\subsection{Exchangeably Weighted Bootstrap}\label{sec:boo}
Define the bootstrap M-estimator
$(\widehat{\theta}^{\ast},\widehat{\eta}^{\ast})
=\arg\sup_{\theta\in\Theta,\eta\in\mathcal{H}}\sum_{i=1}^{n}
m(\theta,\eta)(X_{i}^{\ast})$, where $(X_1^\ast,\ldots,X_n^\ast)$
is the bootstrap sample. Note that the Efron's  nonparametric
bootstrap consists of independent draws with replacement from the
original observations. In this case, we can re-express
$$(\widehat{\theta}^{\ast},\widehat{\eta}^{\ast})
=\arg\sup_{\theta\in\Theta,\eta\in\mathcal{H}}\sum_{i=1}^{n}
W_{ni}m(\theta,\eta)(X_{i}),$$ where
$(W_{n1},\ldots,W_{nn})\sim\mbox{Mult}_n(n,
(n^{-1},\ldots,n^{-1}))$. This multinomial formulation can be
naturally generalized to a class of exchangeable bootstrap weights
$\{W_{ni}\}_{i=1}^n$ whose distribution corresponds to different
bootstrap sampling schemes. This general bootstrap method, called
exchangeably weighted bootstrap, was first proposed by Rubin
(1981) and then extensively studied in \citep{bb95, pw93,mn92}.
The class of exchangeably weighted bootstrap is practically
useful. For example, in Cox regression model, the nonparametric
bootstrap often gives many ties when it is applied to censored
survival data due to its ``discreteness" while the general
weighting scheme comes to the rescue. Other variations of
nonparametric bootstrap are also studied in \citep{cb05} using the
term ``generalized bootstrap".

Let
$\|W_{n1}\|_{2,1}=\int_{0}^{\infty}\sqrt{P_W(|W_{n1}|\geq u)}du$, where $P_W$ is the weight distribution. The bootstrap weights $W_{ni}$'s are assumed to satisfy the
following conditions given in \cite{pw93}:
\begin{enumerate}
\item[W1.] The vector $W_{n}=(W_{n1},\ldots,W_{nn})'$ is
exchangeable for all $n=1,2,\ldots$, i.e., for any permutation
$\pi=(\pi_{1},\ldots,\pi_{n})$ of $(1,2,\ldots,n)$, the joint
distribution of $\pi(W_{n})=(W_{n\pi_{1}},\ldots,W_{n\pi_{n}})'$
is the same as that of $W_{n}$.

\item[W2.] $W_{ni}\geq 0$ for all $n$, $i$ and
$\sum_{i=1}^{n}W_{ni}=n$ for all $n$.

\item[W3.] Assume $\lim\sup_{n\rightarrow\infty}\|W_{n1}\|_{2,1}<\infty$.

\item[W4.]
$\lim_{\lambda\rightarrow\infty}\lim\sup_{n\rightarrow\infty}
\sup_{t\geq\lambda}t^{2}P_W(W_{n1}>t)=0$.

\item[W5.]
$(1/n)\sum_{i=1}^{n}(W_{ni}-1)^{2}\overset{P_W}{\longrightarrow}c^{2}>0$.
\end{enumerate}
Condition W3 is slightly stronger than the bounded second moment
but is implied whenever a $(2+\epsilon)$ absolute moment exists for
some $\epsilon>0$; see Appendix~\ref{ineapp}. By the Markov's
inequality, Condition W4 is satisfied if the $(2+\epsilon')$
moment of $W_{n1}$ is finite for some $\epsilon'>0$. The value of
$c$ depends on the resampling
method, e.g., $c=1$ for nonparametric bootstrap. The bootstrap
weights corresponding to nonparametric bootstrap satisfy W1--W5.
Below, we present several bootstrap examples satisfying W1
-- W5 as shown in Praestgaard and Wellner (1993) where we can find
more details on the sampling schemes.

{\it Example 1. i.i.d.-Weighted Bootstraps}

In this example, the bootstrap weights are defined as
$W_{ni}=\omega_{i}/\bar\omega_n$, where $\omega_{1}, \omega_{2},
\ldots, \omega_{n}$ are i.i.d. positive r.v.s. with
$\|\omega_1\|_{2,1}<\infty$ and $\bar\omega_n=\sum_{i=1}^{n}\omega_i/n$. Thus, we can choose $\omega_i\sim
\mbox{Exponential}(1)$ or $\omega_i\sim\mbox{Gamma}(4, 1)$. The
former corresponds to the Bayesian bootstrap. The multiplier
bootstrap is often thought to be a smooth alternative to the
nonparametric bootstrap; see \citep{l93}. The value of $c^2$ is
calculated as $Var(\omega_1)/(E\omega_1)^2$.

{\it Example 2. The delete-$h$ Jackknife}

In the delete-$h$ jackknife \citep{w87}, the bootstrap weights are
generated by permuting the deterministic weights
\begin{equation}
w_{ni}=\left\{ \begin{array}{cc} \frac{n}{n-h}& 1\leq i\leq n-h,\\
                                   0 & \mbox{otherwise},
                 \end{array}\right.
\end{equation}x
with $\sum_{i=1}^{n}w_{ni}=n$. Specifically, we have $W_{nj}=w_{nR_n(j)}$ where $R_n(\cdot)$ is a
random permutation uniformly distributed over $\{1,\ldots,n\}$. In
Condition W5, $c^2=h/(n-h)$. Thus, we need to choose
$h/n\rightarrow\alpha\in(0,1)$ for $c$ to be positive. Therefore,
Condition W5 does not hold for the ordinary jackknife with $h=1$.

{\it Example 3. The Double Bootstrap}

In the double bootstrap, the bootstrap weights have the following
distribution
\begin{eqnarray}
(W_{n1},\ldots,W_{nn})\sim\mbox{Mult}_n\left(n,(\widetilde
W_{n1}/n,\ldots, \widetilde W_{nn}/n)\right),\label{eg3}
\end{eqnarray}
 conditional on
$\widetilde W_n$ following $\mbox{Mult}_n(n,
(n^{-1,},\ldots,n^{-1}))$. The value of $c$ is $\sqrt{2}$ in this
example.

{\it Example 4. The Polya-Eggenberger Bootstrap}

In this example, the bootstrap weights follow the multinomial
distribution
\begin{eqnarray}
(W_{n1},\ldots,W_{nn})\sim\mbox{Mult}_n\left(n,(D_{n1},\ldots,D_{nn})\right),\label{eg4}
\end{eqnarray}
conditional on
$(D_{n1},\ldots,D_{nn})\sim\mbox{Dirichlet}_n(\alpha,\ldots,\alpha)$
with $\alpha>0$. The value of $c^2$ is calculated as
$(\alpha+1)/\alpha$.

{\it Example 5. The Multivariate Hypergeometric Bootstrap}

As a particular urn-based bootstrap, the bootstrap weights follow
the multivariate hypergeometric distribution with density
\begin{eqnarray}
P(W_{n1}=w_1,\ldots,W_{nn}=w_n)=\frac{\binom{K}{w_1}\cdots\binom{K}{w_n}}{\binom{nK}{n}}\label{eg5}
\end{eqnarray}
for some positive integer $K$. Condition W5 is satisfied with
$c^2=(K-1)/K$.

Under Conditions W1 -- W5 and other regularity conditions, Cheng
and Huang (2010) prove
\begin{eqnarray}
&&(\sqrt{n}/c)(\widehat\theta^\ast-\widehat\theta)\overset{d}{\Longrightarrow}
N(0, \Sigma)\;\;\;\;\;\;\mbox{conditional on $\mathcal{X}_n\equiv
(X_1,\ldots,X_n)$},\label{b-mestdis}
\end{eqnarray}
\noindent where $``\overset{d}{\Longrightarrow}"$ represents the
conditional weak convergence (in probability) defined in
\citep{h84, vw96} (also see (\ref{inter2})) and $P_{W|\mathcal{X}_n}$ is the conditional
probability given $\mathcal{X}_n$. In view of (\ref{b-mestdis}),
the bootstrap variance estimate for $\theta$ is constructed as
\begin{eqnarray}
\widehat\Sigma^\ast=(n/c^2)E_{W|\mathcal{X}_n}(\widehat\theta^\ast-\widehat\theta)(\widehat\theta^\ast-\widehat\theta)',\label{boovar}
\end{eqnarray}
where $E_{W|\mathcal{X}_n}$ is the conditional expectation given
the observed data $\mathcal{X}_n$. We say that the bootstrap
variance estimate is consistent if
$\widehat\Sigma^\ast\overset{P_X}{\longrightarrow}\Sigma$. In
practice, $\widehat\Sigma^\ast$ can be well approximated as
follows:
$${\textstyle
\widehat\Sigma^\ast\approx\widetilde\Sigma^\ast\equiv(n/Bc^2)\sum_{b=1}^{B}\left(\widehat\theta^\ast(b)-B^{-1}
\sum_{b=1}^{B}\widehat\theta^{\ast}(b)\right)\left(\widehat\theta^\ast(b)-B^{-1}
\sum_{b=1}^{B}\widehat\theta^{\ast}(b)\right)'},$$ where
$\widehat\theta^\ast(b)$ is computed based on the $b$-th bootstrap
sample, for sufficiently large number $B$ of bootstrap
repetitions.

\section{Main Result: Bootstrap Moment
Consistency}\label{bootcons} In this section, we will establish
the bootstrap moment consistency of $\theta$ which directly
implies the consistency of $\widehat\Sigma^\ast$ and $t$-type
bootstrap confidence set. To obtain the $p$-th moment consistency
comparing to the distribution consistency, the only additional
cost is to strengthen the $L_1$ maximal inequality condition
required in the latter to the $L_{p'}$ maximal inequality
condition for $p'>1$, i.e., Condition M2. A simple sufficient
condition for M2 i.e., (\ref{boowei}), is also given in terms of
the bootstrap weights, and is verified in the above bootstrap
examples.

It is well known that the convergence in distribution implies the
convergence in moment under the uniform integrability condition.
Lemma 2.1 of Kato (2011) further shows that the above argument is
also valid for the conditional weak convergence in the case of
nonparametric bootstrap. In fact, his arguments (after minor
modifications) can also be applied to the above class of
exchangeably weighted bootstrap; see Lemma~\ref{keylem} below.

\begin{lemma}\label{keylem}
Let $T_n^\ast$ be a scalar statistic of $(X_1,\ldots,X_n)$ and
$(W_{n1},\ldots,W_{nn})$. Suppose that bootstrap weight $W_n$
satisfies W1 -- W5 and the conditional distribution of $T_n^\ast$
given $\mathcal{X}_n$ converges weakly to some fixed distribution
$\mu$ in $P_X$-probability. If
$E_{W|\mathcal{X}_n}|T_n^\ast|^{q'}=O_{P_X}(1)$ for some $q'>1$,
then
$E_{W|\mathcal{X}_n}(T_n^\ast)^{q}\overset{P_X}{\longrightarrow}\int
t^{q} d\mu(t)$ for any integer $1\leq q<q'$.
\end{lemma}

Let $``\aplt"$ $(``\apgt")$ denote smaller (greater) than, up to
an universal constant. Denote $E_{XW}$ and $P_{XW}$ as the joint
expectation and joint probability, respectively. Let $P_X f=\int
fdP_X$, $\mathbb{P}_nf=\sum_{i=1}^n f(X_i)/n$ and
$\mathbb{P}_n^\ast f=\sum_{i=1}^n f(X_i^\ast)/n=\sum_{i=1}^n
W_{ni}f(X_i)/n$. For example, we can rewrite
$(\widehat\theta,\widehat\eta)=\arg\sup_{\theta\in\Theta,\eta\in\mathcal{H}}\mathbb{P}_n
m(\theta,\eta)$ and
$(\widehat\theta^\ast,\widehat\eta^\ast)=\arg\sup_{\theta\in\Theta,\eta\in\mathcal{H}}\mathbb{P}_n^\ast
m(\theta,\eta)$. Define the empirical process $$\mathbb{G}_nf=\sqrt{n}(\mathbb{P}_n-P_X)f$$ and its norm $$\|\mathbb{G}_n\|_{\mathcal{F}}=\sup_{f\in\mathcal{F}}|\mathbb{G}_nf|$$ as well as their bootstrapped analogues $$\mathbb{G}_n^\ast f=\sqrt{n}(\mathbb{P}_n^\ast-\mathbb{P}_n)f$$ and
$$\|\mathbb{G}_n^\ast\|_{\mathcal{F}}=
\sup_{f\in\mathcal{F}}|\mathbb{G}_n^\ast f|.$$ For
any class of functions $\mathcal{A}$ under a metric $\ell$, we
define $\log N_{[]}(\epsilon, \mathcal{A},\ell)$ and $\log
N(\epsilon, \mathcal{A}, \ell)$ as the $\epsilon$-bracketing entropy
number and $\epsilon$-entropy number, respectively. The related
bracketing entropy integral and uniform entropy integral are thus
\begin{eqnarray*}
J_{[]}(\delta,\mathcal{A},\ell)&=&\int_{0}^{\delta}\sqrt{1+\log
N_{[]}(\epsilon,\mathcal{A},\ell)}d\epsilon,\\
J(\delta,\mathcal{A})&=&\sup_{Q}\int_0^\delta\sqrt{1+\log
N(\epsilon\|A\|_{L_2(Q)}, \mathcal{A}, L_2(Q))}d\epsilon,
\end{eqnarray*}
where $A$ is the envelope function of $\mathcal{A}$, and the
supreme is taken over all discrete probability measures $Q$ with
$\|A\|_{L_2(Q)}>0$.

In the following, we provide a set of sufficient conditions for
bootstrap moment consistency.
\begin{enumerate}
\item[M1.] For any $(\theta,\eta)\in\Theta\times\mathcal{H}$, we
have
\begin{eqnarray}
E_X(m(\theta,\eta)-m(\theta_0,\eta_{0}))\aplt-\|\theta-\theta_{0}\|^{2}-d^{2}
(\eta,\eta_{0}).\label{m1con1}
\end{eqnarray}


\item[M2.] Define
$\mathcal{N}_\delta=\{m(\theta,\eta)-m(\theta_0,\eta_0):
\|\theta-\theta_0\|\leq\delta, d(\eta,\eta_0)\leq\delta,
(\theta,\eta)\in\Theta\times\mathcal{H}\}$. We assume that, for
some $p'\geq 1$ and every $\delta>0$,
\begin{eqnarray}
\left(E_X\|\mathbb{G}_n\|_{\mathcal{N}_\delta}^{p'}\right)^{1/p'}&\aplt&\delta\label{conmod3},\\
\left(E_{XW}\|\mathbb{G}_n^\ast\|_{\mathcal{N}_\delta}^{p'}\right)^{1/p'}&\aplt&\delta.\label{conmod4}
\end{eqnarray}

\item[M3.] Assume that $d(\widehat\eta,\eta_0)=O_{P_X}(n^{-1/2})$
and $d(\widehat\eta^\ast,\eta_0)=O_{P_{XW}}(n^{-1/2})$.
\end{enumerate}

Let $\rho_p(\cdot): \mathbb R^d\mapsto\mathbb R$ be a continuous function with polynomial growth rate $p$, i.e., there exist constants $K_1, K_2>0$ and $p>1$ such that, for any $x$, $$|\rho_p(x)|\leq K_1+K_2\|x\|^p.$$
\begin{theorem}\label{bconthm}
Suppose that Conditions W1 -- W5 and M1 -- M3 hold. If
$\widehat\theta^\ast$ is distribution consistent, i.e.,
(\ref{b-mestdis}), then we have
\begin{eqnarray}
E_{W|\mathcal{X}_n}\rho_p(\sqrt{n}(\widehat\theta^\ast-\widehat\theta))\overset{P_X}
{\longrightarrow}E\rho_p(T),\label{mainres}
\end{eqnarray}
where $T\sim N(0,\Sigma)$, for any integer $p$ satisfying $1\leq
p<p'$.
\end{theorem}
\noindent Note that the consistency of bootstrap variance estimate is immediately implied by the above Theorem when $p=2$.

Condition M1 assumes the quadratic behavior of the criterion
function $(\theta,\eta)\mapsto E_X m(\theta, \eta)$. Condition M2
assumes two maximal inequalities in terms of $L_{p'}$-norm for
$p'\geq 1$. Both conditions impose global restrictions on the criterion function. In comparison, $\sqrt{n}$-consistency Theorem in Page 291 of \cite{vw96} only requires their local counterparts. The global type conditions are absolutely needed for obtaining the moment consistency since we need to control the behavior of $P_{XW}(\sqrt{n}(\widehat\theta^\ast-\theta_0)>t)$ for large $t$ and also control the behavior of $\mathbb G_n^\ast(m(\theta,\eta)-m(\theta_0,\eta_0))$ over all ``shells" $S_{j,n}$ (defined in (\ref{shell})) partitioning $\Theta\times\mathcal H$. The
convergence rate of the bootstrap estimate in Condition M3, i.e.,
$d(\widehat\eta^\ast,\eta_0)=O_{P_{XW}}(n^{-1/2})$, can also be
understood in the following way: for any $\delta>0$, there exists
a $0<L<\infty$ such that
$$P_X\left(P_{W|\mathcal{X}_n}\left(\sqrt{n}d(\widehat\eta^\ast,\eta_0)\geq L\right)>\delta\right)\longrightarrow 0\;\;\;\;\;
\mbox{as}\;n\rightarrow\infty.$$ We can verify Condition M3 using
Theorem 2 of \citep{ch10} under very weak model assumptions. On the other hand, we have to admit that relaxing root-n rate requirement is quite challenging. The key technical reason is that the shelling argument in the proof of Theorem~\ref{bconthm} requires partitioning the parameter space $\Theta\times\mathcal H$ as a whole (and thus both estimators are required to have the same root-n rate of convergence); see the shell definition in (\ref{shell}). The separate partition of $\Theta$ and $\mathcal H$ seems to be a solution, but this idea cannot be easily adopted in the current framework of proof.


Below, we discuss three different approaches for verifying
(\ref{conmod3}). Lemma 2.14.1 in \citep{vw96} implies that
\begin{eqnarray}
\left(E_X\|\mathbb{G}_n\|_{\mathcal{N}_\delta}^{p'}\right)^{1/p'}\aplt
J(1,\mathcal{N}_\delta)\|N_\delta\|_{L_{2\vee p'}(P_X)},
\end{eqnarray}
where $N_\delta$ is the envelope function of $\mathcal{N}_\delta$.
Thus, Condition (\ref{conmod3}) holds if
\begin{eqnarray}
J(1,\mathcal{N}_\delta)&<&\infty,\label{entcon}\\
\|N_\delta\|_{L_{2\vee p'}(P_X)}&\aplt&\delta.\label{int2}
\end{eqnarray}
The typical function classes with finite uniform entropy integral
include the VC class and the related larger VC-hull class; see
their definitions in Section 2.6 of \citep{vw96}. Under the
(global) Lipschitz continuous condition:
\begin{eqnarray}
|m(\theta,\eta)(x)-m(\theta_0,\eta_0)(x)|\leq
M(x)(\|\theta-\theta_0\|+d(\eta,\eta_0))\;\;\;\mbox{for
any}\;(\theta,\eta)\in\Theta\times\mathcal{H},\label{m2ver5}
\end{eqnarray}
we can show (\ref{int2}) if $E_XM^{2\vee p'}(X)<\infty$. The above global Lipschitz condition (\ref{m2ver5}) (together with M1) is usually easy to verify given that $(\theta_0, \eta_0)$ is true value. Alternatively, by decomposing
$(m(\theta,\eta)-m(\theta_0,\eta_0))$ as the sum of
$(m(\theta,\eta)-m(\theta_0,\eta))$ and
$(m(\theta_0,\eta)-m(\theta_0,\eta_0))$, we can also verify
(\ref{conmod3}) if the following holds:
$$\left(E_X\|\mathbb{G}_n\|_{\dot{\mathcal{N}}_{\delta}}^{p'}\right)^{1/p'}<\infty\;\;\;\mbox{and}\;\;\;\left(E_X\|\mathbb{G}_n\|_{\mathcal{N}_{\delta2}}^{p'}\right)^{1/p'}\aplt\delta,$$
where
$\dot{\mathcal{N}}_\delta=\{(\partial/\partial\theta)m(\theta,\eta):
\|\theta-\theta_0\|\leq\delta, d(\eta,\eta_0)\leq\delta \}$ and
$\mathcal{N}_{\delta2}=\{m(\theta_0,\eta)-m(\theta_0,\eta_0):
d(\eta,\eta_0)\leq\delta\}$. Again, Lemma 2.14.1 in \citep{vw96}
can be applied here. Our third approach is to bound the higher
moments $(E_X\|\mathbb{G}_n\|_{\mathcal{N}_\delta}^{p'})^{1/p'}$
for $p'>1$ by $E_X\|\mathbb{G}_n\|_{\mathcal{N}_\delta}$ plus some
norm of $N_\delta$, based on the following two inequalities:
\begin{eqnarray}
\left(E_X\|\mathbb{G}_n\|_{\mathcal{N}_\delta}^{p'}\right)^{1/p'}&\aplt&
E_X\|\mathbb{G}_n\|_{\mathcal{N}_\delta}+n^{1/2-1/p'}\|N_{\delta}\|_{\psi_{p'}}\;\;\;\mbox{for}\;1<p'<2,\label{m2ver}\\
\left(E_X\|\mathbb{G}_n\|_{\mathcal{N}_\delta}^{p'}\right)^{1/p'}&\aplt&
E_X\|\mathbb{G}_n\|_{\mathcal{N}_\delta}+n^{-1/2+1/p'}\left(E_X|N_\delta|^{p'}\right)^{1/p'}\;\;\;\mbox{for}\;p'\geq
2,\label{m2ver2}
\end{eqnarray}
where $\|\cdot\|_{\psi_p}$ is the Orlicz norm with
$\psi_p(t)=\exp(t^p)-1$. The above two inequalities are derived
based on Theorem 2.14.5 in \citep{vw96} and the fact that the
$\psi_p$-norm dominates the $L_p$-norm for each $p$. Now, we
assume (\ref{m2ver5}). When $p'>1$ but $\neq 2$, the second term
in the right hand side of (\ref{m2ver}) ((\ref{m2ver2})) converges
to zero as $n\rightarrow\infty$ if $\|M\|_{\psi_{p'}}<\infty$
$(\|M\|_{L_{p'}(P_X)}<\infty)$. When $p'=2$, the second term in
the right hand side of (\ref{m2ver2}) is of the order $O(\delta)$
if $\|M\|_{L_{2}(P_X)}<\infty$. Thus, if
$E_X\|\mathbb{G}_n\|_{\mathcal{N}_\delta}\aplt\delta$, we can show
(\ref{conmod3}). Fortunately, several technical tools are
available to compute the upper bound of
$E_X\|\mathbb{G}_n\|_{\mathcal{N}_\delta}$ in terms of the
bracketing entropy integral (using Theorem 2.14.2 or Lemma 3.4.2
in \citep{vw96}) or the uniform entropy integral (see van der
Vaart and Wellner (2011)). For example, in view of the above
analysis and Theorem 2.14.2 in \citep{vw96}, a simple sufficient
condition for (\ref{conmod3}) is
$$J_{[]}(1, \mathcal{N}_\delta, L_2(P_X))+\|M\|_{\psi_{p'\vee 2}}<\infty\;\mbox{and Condition}\;(\ref{m2ver5})$$ due to the fact that the
$\psi_p$-norm dominates the $L_p$-norm for each $p$, and
$\psi_q$-norm for any $q\leq p$.

To verify (\ref{conmod4}), we will employ the general $L_p$
multiplier inequality developed in Appendix~\ref{lpmin} to bound
$(E_{XW}\|\mathbb{G}_n^\ast\|_{\mathcal{N}_\delta}^{p'})^{1/p'}$.
According to Appendix~\ref{conver}, it suffices to show the
following bootstrap weight condition
\begin{eqnarray}
W_{n1}^{p'}\;\;\mbox{satisfies}\;\;\mbox{Conditions}\;\;W3\; \&\;
W4;\label{boowei}
\end{eqnarray}
if (\ref{conmod3}) holds. Condition (\ref{boowei}) is essentially
very weak; see discussions in {\em Examples 1 -- 5} below. In the
end, we want to point out that Conditions W1 -- W5 and M1 -- M3
(when $p'=1$) are also needed in showing the bootstrap
distribution consistency (\ref{b-mestdis}); see Theorems 1 \& 3 of
\citep{ch10}. In view of the above discussions, it appears that we
only need to strengthen the $L_1$ maximal inequalities to the
$L_{p'}$ maximal inequalities for $p'\geq 1$ to achieve the
bootstrap moment consistency beyond the distribution consistency.

An obvious implication of Theorem~\ref{bconthm} is that the
bootstrap moment estimate of arbitrary order is consistent if
Condition M2 is valid for all $p'\geq 1$. It is worthwhile to
remark that the uniform integrability of $\widehat\theta$, i.e.,
$E_X\|\sqrt{n}(\widehat\theta-\theta_0)\|^{p}<\infty$, is also
proven in the proof of Theorem~\ref{bconthm}. Thus, under the same
set of conditions, the moment convergence of $\widehat\theta$ also
follows. In addition, Theorem~\ref{bconthm} is also valid even for
the approximate maximizer, i.e.,
\begin{eqnarray*}
\mathbb{P}_{n}m(\widehat\theta,\widehat{\eta})&\geq&
\mathbb{P}_{n}m(\theta_0,\eta_{0})-O_{P_{X}}(n^{-1}),\\
\mathbb{P}_{n}^{\ast}m(\widehat\theta^\ast,\widehat{\eta}^{\ast})&\geq&
\mathbb{P}_{n}^{\ast}m(\theta_0,\eta_{0})-O_{P_{XW}}(n^{-1})
\end{eqnarray*}
after slightly modifying its proof.

The distribution consistency result (\ref{b-mestdis}) directly
implies the consistency of bootstrap hybrid and percentile
confidence sets. Given the consistent variance estimate
$\widehat{\Sigma}$ based on $(X_1,\ldots,X_n)$, the more accurate
{\it t}-type bootstrap confidence set is constructed as
$$BC_{t}(\alpha)=\left[\widehat{\theta}-\frac{\widehat{\Sigma}^{1/2}
\omega_{n(1-\alpha/2)} ^{\ast}}{\sqrt{n}},\widehat{\theta}-\frac{
\widehat{\Sigma}^{1/2}
\omega_{n(\alpha/2)}^{\ast}}{\sqrt{n}}\right],$$ where $\omega
_{n\alpha}^{\ast}$ satisfies
$P_{W|\mathcal{X}_n}((\sqrt{n}/c)(\widehat{\Sigma}^{\ast})^{-1/2}(\widehat{\theta}
^{\ast}-\widehat{\theta})\leq\omega _{n\alpha}^{\ast})=\alpha$ and $``\leq"$ is componentwise.
Note that $\omega_{n\alpha}^\ast$ is not unique when $\theta$ is a
vector.  The following Corollary theoretically justifies the
widely used bootstrap variance estimate $\widehat\Sigma^\ast$, and
further establishes the consistency of $t$-type confidence set
$BC_t(\alpha)$.
\begin{corollary}\label{coro}
Suppose that Conditions in Theorem~\ref{bconthm} hold. If we
further require that Condition M2 holds for some $p'>2$, then we
have
\begin{eqnarray}
\widehat\Sigma^\ast&\overset{P_X}{\longrightarrow}&\Sigma,\label{coro1}\\
P_{XW}(\theta_0\in BC_{t}(\alpha))&\longrightarrow&
1-\alpha\label{coro2}
\end{eqnarray}
as $n\rightarrow\infty$.
\end{corollary}
The variance consistency (\ref{coro1}) directly follows from
Theorem~\ref{bconthm}. To show the consistency of $t$-type
confidence set, i.e., (\ref{coro2}), we apply the Slutsky's Lemma
and its conditional version given in Appendix~\ref{conslt}
(together with Lemma 4.6 of \citep{pw93}) to (\ref{mestdis}) and
(\ref{b-mestdis}). Thus, for any fixed $x\in\mathbb{R}^d$, we
obtain that
\begin{eqnarray}
P_X(\sqrt{n}\widehat\Sigma^{-1/2}(\widehat{\theta}-\theta_0)\leq
x)&\longrightarrow&
\Psi(x),\label{ab}\\
P_{W|\mathcal{X}_n}((\sqrt{n}/c)(\widehat\Sigma^\ast)^{-1/2}(\widehat{\theta}^{\ast}-
\widehat{\theta})\leq x)&\overset{P_X}
{\longrightarrow}&\Psi(x),\label{ab2}
\end{eqnarray}
where $\Psi(x)=P(N(0,I)\leq x)$. A straightforward application of
Lemma 23.3 in \citep{v98} concludes the proof of (\ref{coro2})
based on (\ref{ab}) \& (\ref{ab2}).

In the end of this section, we will verify the bootstrap weight
condition (\ref{boowei}) in six different types of bootstraps
introduced in Section~\ref{sec:boo}.

{\it Example 1. i.i.d.-Weighted Bootstraps (Cont')}

We will show that (\ref{boowei}) holds under the assumption that
$\omega_i$ has bounded $(2+\epsilon)p'$-th moment for some
$\epsilon>0$. This assumption implies that
\begin{eqnarray}
\|\omega_i^{p'}\|_{2,1}<\infty\label{int1}
\end{eqnarray}
based on Appendix~\ref{ineapp}. The derivations in Page 2080 of
\citep{pw93} give that
\begin{eqnarray}
P_W(W_{n1}^{p'}>t)\leq
P(\omega_1>t^{1/p'}(1-\epsilon))+t^{-p/(2p')}n^{p/2}\rho(\epsilon)^{n/2}\label{inter6}
\end{eqnarray}
for any $0<\epsilon<1$, $p>0$ and some $0<\rho(\epsilon)<1$, which
further implies that
\begin{eqnarray*}
\|W_{n1}^{p'}\|_{2,1}&\leq&
1+\int_1^\infty\sqrt{P_W(W_{n1}^{p'}>t)}dt\\
&\leq&1+\int_1^{\infty}\sqrt{P(\omega_1^{p'}>t(1-\epsilon)^{p'})}dt+n^{p/4}\rho(\epsilon)^{n/4}\int_1^\infty
t^{-p/4p'}dt\\
&\leq&1+\frac{1}{(1-\epsilon)^{p'}}\|\omega_1^{p'}\|_{2,1}+n^{p/4}\rho(\epsilon)^{n/4}\int_1^\infty
t^{-p/4p'}dt.
\end{eqnarray*}
By choosing $p>4p'$, we know that
$\lim\sup_{n\rightarrow\infty}\|W_{n1}^{p'}\|_{2,1}<\infty$ due to
(\ref{int1}). To see that $W_{n1}^{p'}$ satisfies Condition W4, it
suffices to show that $\lim_{t\rightarrow\infty}
t^2P(\omega_1^{p'}>t)=0$ according to (\ref{inter6}). This is
implied by the Markov's inequality and the bounded moment
assumption on $\epsilon$.

{\it Example 2. The delete-$h$ Jackknife (Cont')}

Recall that the bootstrap weight $W_{nj}=w_{nR_n(j)}$. Then, we
have
\begin{eqnarray}
P_W(W_{n1}>t)=\sharp\{j:
w_{nj}>t\}=\frac{n-h}{n}1\{t<n/(n-h)\}.\label{eg2-0}
\end{eqnarray}
In view of (\ref{eg2-0}), Condition (\ref{boowei}) can be verified
as follows
\begin{eqnarray}
&&\lim\sup_{n\rightarrow\infty}\int_0^\infty
\sqrt{P_W(W_{n1}>u)}du^{p'}=\lim\sup_{n\rightarrow\infty}
\left(\frac{n}{n-h}\right)^{p'-1/2}=\left(\frac{1}{1-\alpha}\right)^{p'-1/2}<\infty,\label{eg2-1}\\
&&\lim_{t\rightarrow\infty}\lim\sup_{n\rightarrow\infty}\frac{n-h}{n}t^21\left\{t<\frac{n^{p'}}{(n-h)^{p'}}\right\}=\lim_{t\rightarrow\infty}(1-\alpha)t^2
1\{t<(1-\alpha)^{-p'}\}=0.\label{eg2-2}
\end{eqnarray}

A sufficient condition for (\ref{boowei}) is
\begin{eqnarray}
\lim\sup_{n\rightarrow\infty}E_WW_{n1}^{(2+\epsilon)p'}<\infty\;\;\;\mbox{for
some}\;\epsilon>0.\label{boowei2}
\end{eqnarray}
This can be proven based on the Appendix~\ref{ineapp} and
Chebyshev's inequality as remarked above. Thus, to guarantee the
bootstrap variance consistency, i.e. Corollary~\ref{coro}, we only
need to require
\begin{eqnarray}
\lim\sup_{n\rightarrow\infty}E_WW_{n1}^{5}<\infty\label{boowei3}
\end{eqnarray}
since we can always choose $p'=5/(2+\epsilon)>2$ for small enough
$\epsilon>0$. Assuming
$W_n=(W_{n1},\ldots,W_{nn})'=\mbox{Mult}_n(n, (p_1,\ldots,p_n))$,
we have
\begin{eqnarray}
E_WW_{n1}^5=np_1+15n^{(2)}p_1^2+25n^{(3)}p_1^3+10n^{(4)}p_1^4+n^{(5)}p_1^5,\label{mult}
\end{eqnarray}
where $n^{(k)}=n(n-1)\cdots(n-k+1)$, according to Page 33 in
\cite{jkb97}. If $p_i=1/n$ for $i=1,\ldots,n$, we know
$E_WW_{n1}^5<52$. Thus, Condition (\ref{boowei3}) (also
(\ref{boowei})) is trivially satisfied in the {\em Efron's
nonparametric bootstrap}. Condition (\ref{boowei3}) can be easily
verified in the examples 3 -- 5 discussed before.

{\it Example 3. The Double Bootstrap (Cont')}

Based on (\ref{eg3}) \& (\ref{mult}), we can compute $E_WW_{n1}^5$
as
\begin{eqnarray*}
&&E(E_W(W_{n1}^5|\widetilde W_n))\\&=&E\left(\widetilde
W_{n1}+15(n^{(2)}/n^2)\widetilde
W_{n1}^2+25(n^{(3)}/n^3)\widetilde
W_{n1}^3+10(n^{(4)}/n^4)\widetilde
W_{n1}^{4}+(n^{(5)}/n^5)\widetilde W_{n1}^5\right),
\end{eqnarray*}
which implies Condition (\ref{boowei3}) since $E\widetilde
W_{n1}^5<52$.

{\it Example 4. The Polya-Eggenberger Bootstrap (Cont')}

Following similar analysis in double bootstrap and (\ref{eg4}), we
have
$$E_WW_{n1}^5=E\left(nD_{n1}+15n^{(2)}D_{n1}^2+25n^{(3)}D_{n1}^3+10n^{(4)}D_{n1}^4+n^{(5)}D_{n1}^5\right).$$
We can verify (\ref{boowei3}) if we can show
$$\lim\sup_{n\rightarrow\infty}n^{(p)}ED_{n1}^p<\infty$$ for
$p=1,\ldots,5$. This is essentially true for all $p$ based on the
below derivations
\begin{eqnarray*}
n^{(p)}ED_{n1}^p=n^{(p)}\frac{\alpha\cdots(\alpha+p-1)}{n\alpha\cdots(n\alpha+p-1)}\longrightarrow\prod_{k=1}^{p-1}\frac{\alpha+k}{\alpha}\;\;\;\mbox{as}\;n\rightarrow\infty,
\end{eqnarray*}
where the formula for calculating $ED_{n1}^p$ is given in Page 96
of \citep{jk77}.

{\it Example 5. The Multivariate Hypergeometric Bootstrap (Cont')}

According to (\ref{eg5}) and Page 96 of \citep{jk77}, we have
$$E_{W}W_{n1}^5=a_{n,K}(1)+15a_{n,K}(2)+25a_{n,K}(3)+10a_{n,K}(4)+a_{n,K}(5),$$ where
$a_{n,K}(r)=n^{(r)}K^{(r)}/(nK)^{(r)}$. Since
$a_{n,K}(r)<K^{(r)}$, we can show
$\lim\sup_{n\rightarrow\infty}E_WW_{n1}^5<\infty$.

%
%
%
%
%
%
%
%
%
%
%
%
%
%
%
%
%

\section{Cox Regression Model with Right Censored Data}\label{example}
We use the following Cox regression model to illustrate the practicality of the stated conditions M1 -- M3, and then run simulations for the five classes of bootstrap methods considered in Examples 1 -- 5. Indeed, the advantages of using bootstrap inferences in this model were
considered in the literature, e.g., \citep{et86}. In the Cox
regression model, the hazard function of the survival time $T$ of
a subject with covariate $Z$ is modelled as:
\begin{eqnarray}
\lambda(t|z)\equiv\lim_{\Delta\rightarrow
0}\frac{1}{\Delta}P(t\leq T<t+\Delta|T\geq
t,Z=z)=\lambda(t)\exp(\theta' z),\label{eg1den}
\end{eqnarray}
where $\lambda$ is an unspecified baseline hazard function and
$\theta$ is a regression vector. In this model, we are usually
interested in $\theta$ while treating the cumulative hazard
function $\eta(y)=\int_{0}^{y}\lambda(t)dt$ as the nuisance
parameter. With right censoring of survival time, the data
observed is $X=(Y,\delta,Z)$, where $Y=T\wedge C$, $C$ is a
censoring time, $\delta=I\{T\leq C\}$, and $Z$ is a regression
covariate belonging to a compact set
$\mathbb{C}\subset\mathbb{R}^d$. We assume that $C$ is independent
of $T$ given $Z$. The log-likelihood is thus
\begin{eqnarray}
m(\theta,\eta)(x)=\delta\theta' z-\exp(\theta'
z)\eta(y)+\delta\log\eta\{y\},
\end{eqnarray}
where $\eta\{y\} = \eta(y) - \eta(y-)$ is a point mass that
denotes the jump of $\eta$ at point $y$. The parameter space
$\mathcal{H}$ is restricted to a set of nondecreasing cadlag
functions on the interval $[0,\tau]$ with $\eta(\tau)\leq M$ for
some constant $M$. It is well known that the MLE $\widehat\theta$
is semiparametric efficient with the asymptotic variance obtained
in \citep{bhhw83}:
\begin{eqnarray}
\Sigma=\widetilde I_0^{-1}\equiv
\left\{E\widetilde\ell_{\theta_0,\eta_0}(X)\widetilde\ell_{\theta_0,\eta_0}'(X)\right\}^{-1},
\label{eg1-var}
\end{eqnarray}
where the efficient information matrix $\widetilde I_0$ is
computed via the efficient score function
$$\widetilde\ell_{\theta,\eta}(x)=\delta\left(z-\frac{Ee^{\theta'Z}Z1\{Y\geq y\}}{Ee^{\theta'Z}1\{Y\geq y\}}\right)-
e^{\theta'z}\int_0^y\left(z-\frac{Ee^{\theta'Z}Z1\{Y\geq
t\}}{Ee^{\theta'Z}1\{Y\geq t\}}\right)d\eta(t).$$ The negative
second derivative of the partial likelihood can be used to
estimate $\Sigma^{-1}$. This is a special case of the observed
profile information defined as the negative second numerical
derivative of the profile likelihood; see \citep{mv99}. In
general, this approach requires a careful choice of the step size
and crucially depends on the curvature structure of the profile
likelihood which may not behave well under small sample.

Cheng and Huang (2010) have shown that the exchangeably weighted
bootstrap is consistent in estimating the limiting distribution of
$\widehat\theta$. Below, we will verify that Conditions M1 -- M3
hold for this model such that the bootstrap is also consistent for
estimating $\Sigma$. Since the true value $(\theta_0,\eta_0)$ is
the maximizer of $(\theta,\eta)\mapsto P_Xm(\theta,\eta)$ (under
certain identifiability condition), it is not difficult to verify
Condition M1 by defining $d(\eta,\eta_0)=\|\eta-\eta_0\|_\infty$,
where $\|\cdot\|_{\infty}$ denotes the supreme norm. The
convergence rates of $\widehat\eta$ and $\widehat\eta^\ast$ are
established in Theorem 3.1 of \citep{mv99} and Theorem 2 of
\citep{ch10}, respectively, as
\begin{eqnarray}
\|\widehat{\eta}
-\eta_{0}\|_{\infty}=O_{P_X}(n^{-\frac{1}{2}})\;\;\;\mbox{and}\;\;\;\|\widehat{\eta}^\ast
-\eta_{0}\|_{\infty}=O_{P_{XW}}(n^{-\frac{1}{2}}).\label{eg1crate}
\end{eqnarray}
Thus, we have verified Condition M3. To verify (\ref{conmod3}) in
M2, we apply the first approach by showing (\ref{entcon}) \&
(\ref{int2}). Note that the class of bounded monotone functions,
e.g., $\eta(y)$ and $\eta(y-)$, is VC-hull class. Considering the
form of $m(\theta,\eta)$ (writing $\eta\{y\}=\eta(y)-\eta(y-)$),
we know that (\ref{entcon}) is satisfied by the stability property
of the BUEI (bounded uniform entropy integral) function class, i.e., Lemma 9.14 of \citep{k08}. Note
that (\ref{int2}) trivially holds since we can show (\ref{m2ver5})
with $M(x)$ as some finite constant due to the compactness of
$\mathbb{C}$ and $\mathcal{H}$. This also justifies
$\|N_\delta\|_{L_{p'}(P_X)}<\infty$. Thus, (\ref{conmod4}) holds
according to Appendix~\ref{conver}.

We conclude this section by running simulations for the above five classes of bootstrap methods, and also try to give advice in choosing bootstrap weights accordingly. We consider four different settings, and set $\lambda_0(t)=\exp(t)$ in the simulations. In each setting, $n=500, 1000, 1500$ were generated and the variance was calculated as an average of $100$ replications. The censoring time $C$ follows $U[0,t_n]$ where $t_n$ was chosen such that the average effective sample size over $500, 1000, 1500$ samples is approximates $p\times n$. We applied five different bootstrap methods as specified above: Efron's bootstrap, delete-$h$ Jackknife with $h=0.002n$, double bootstrap, Polya-Eggenberger bootstrap with parameter $\alpha$ and multivariate hypergeometric bootstrap with parameter $K$. The bootstrap variance estimates were calculated based on $1000$ bootstrap repetitions. We used the ``coxph" function in the R package ``survival" to calculate the MLE of $\theta$ and its corresponding variance, which is used as a benchmark. In the setting I, the covariate $z$ was generated from U[0,1], the regression coefficient was set as $\theta_0=2$ and the parameters $p$, $\alpha$ and $K$ were set to be $0.9$, $3$ and $3$. In the setting II, we let $\theta_0=4, p=0.7, \alpha=5$ and $K=4$. In the setting III, we generated $z$ evenly in the interval $[1,2]$ and set $\theta_0=3, p=0.8$, all the other setup are the same as setting II. We consider two dimensional $\theta_0=(2,1)^T$ in the setting IV. The covariate vector $z$ follows two independent uniform distributions: $U[0,1]\times U[0,1]$. We inherited other setup from setting III. All the results are summarized in Table~1.

Given our consistency results, it is not surprising to see that all these exchangeably weighted bootstrap methods produce fairly close results to the variance (covariance matrix) of MLE in all the setup. Their subtle differences are mainly due to the specifications on the data generation mechanism and the choice of bootstrap parameters, e.g., $\alpha$. Having said that, we would like to recommend Efron's bootstrap for practical use since it is the most straightforward to implement with the least computational cost (in contrast with the other four methods). On the other hand, the above observations strongly motivate the second order theoretical studies that may lead to a more refined practical guidance in selecting bootstrap methods. We leave this as a future topic.

\begin{center}
\begin{table}
\caption{Bootstrap variance estimate based on different bootstraps}
\scriptsize{
\begin{tabular}{l|l|c|c|c|c}
\hline\hline
 & & I & II & III & IV\\
\hline\hline
\multirow{12}{*}{n=500} &Maximum Likelihood  & 23.41 & 41.18 & 29.15 & $\left(\begin{array}{cc} 27.51 & 2.37\\  2.37& 24.60\end{array}\right)$\\
\hhline{~-----}
&Efron's bootstrap& 23.22 & \textbf{41.23} & 29.58 & $\left(\begin{array}{cc} 27.88& 2.35\\ 2.35& 24.97\end{array}\right)$\\
\hhline{~-----}
&Delete-$h$ Jackknife & \textbf{23.51} & 41.73& 29.89 & $\left(\begin{array}{cc} 28.01& 2.42\\ 2.42& 25.01\end{array}\right)$\\
\hhline{~-----}
&Double bootstrap & 23.60 & 42.00 & 29.68& $\left(\begin{array}{cc} 28.25& 2.42\\ 2.42& 25.35\end{array}\right)$\\
\hhline{~-----}
&Polya-Eggenberger bootstrap& 23.58 & 41.41 & 29.63 & $\left(\begin{array}{cc} 27.71& 2.23\\ 2.23& 24.78\end{array}\right)$\\
\hhline{~-----}
&Multivariate hypergeometric bootstrap& \textbf{23.31} & 41.42 & \textbf{29.39} &$\left(\begin{array}{cc} \textbf{27.64}&\textbf{ 2.30}\\ \textbf{2.30}& \textbf{24.86}\end{array}\right)$\\
\hline\hline
\multirow{12}{*}{n=1000} &Maximum Likelihood  & 23.41 & 41.66 & 28.72 & $\left(\begin{array}{cc} 26.84& 1.78\\ 1.78& 24.32\end{array}\right)$\\
\hhline{~-----}
&Efron's bootstrap& 23.32 & 41.70 & 29.01 & $\left(\begin{array}{cc}  27.16& 1.72\\ 1.72 &24.61\end{array}\right)$\\
\hhline{~-----}
&Delete-$h$ Jackknife & \textbf{23.37} & 41.79 & 29.07 & $\left(\begin{array}{cc} 27.31& 1.81\\ 1.81 &24.66\end{array}\right)$\\
\hhline{~-----}
&Double bootstrap & 23.33 & 41.49 & 29.06 & $\left(\begin{array}{cc} 27.28& 1.84\\  1.84& 24.89\end{array}\right)$\\
\hhline{~-----}
&Polya-Eggenberger bootstrap& 23.09 & \textbf{41.65} & 28.95 & $\left(\begin{array}{cc} \textbf{27.22}& \textbf{1.80}\\ \textbf{1.80}& \textbf{24.54}\end{array}\right)$\\
\hhline{~-----}
&Multivariate hypergeometric bootstrap& 23.35 & 41.35 & \textbf{28.92} &$\left(\begin{array}{cc} 27.64& 2.30\\ 2.30& 24.86\end{array}\right)$\\
\hline\hline
\multirow{12}{*}{n=1500} &Maximum Likelihood  & 23.20 & 40.81 & 29.15 & $\left(\begin{array}{cc} 26.83& 1.85 \\ 1.85&24.06\end{array}\right)$\\
\hhline{~-----}
&Efron's bootstrap& 23.34 & \textbf{40.75} & \textbf{29.30} & $\left(\begin{array}{cc}  26.79 & 1.94\\ 1.94& 24.13\end{array}\right)$\\
\hhline{~-----}
&Delete-$h$ Jackknife & 23.27 & 40.89 & 29.38 & $\left(\begin{array}{cc} 27.12 & 1.97\\  1.97& 24.25\end{array}\right)$\\
\hhline{~-----}
&Double bootstrap & 23.53 & 41.02 & 29.36 & $\left(\begin{array}{cc} \textbf{27.11} & \textbf{1.84} \\ \textbf{1.84}& \textbf{24.04}\end{array}\right)$\\
\hhline{~-----}
&Polya-Eggenberger bootstrap& 22.92 & 40.62 & 29.31 & $\left(\begin{array}{cc} 27.12& 2.01 \\ 2.01& 24.16\end{array}\right)$\\
\hhline{~-----}
&Multivariate hypergeometric bootstrap& \textbf{23.23} & 40.60 & 29.36 &$\left(\begin{array}{cc} 27.03 & 1.99 \\ 1.99& 24.19\end{array}\right)$\\
\hline\hline
\end{tabular}
}
\end{table}
\end{center}

{\bf Acknowledgment.} The author thanks Professor Yoichi Nishiyama
for sending me his technical note attached to Nishiyama (2010) and
thanks Professor Jon Wellner for helpful discussions.

\address{Department of Statistics\\
Purdue University\\
250 N. University Street\\
West Lafayette, IN 47906\\Email: chengg@purdue.edu}

\section*{Appendix}
\setcounter{subsection}{0}
\renewcommand{\thesubsection}{A.\arabic{subsection}}
\setcounter{equation}{0}
\renewcommand{\theequation}{A.\arabic{equation}}
\setcounter{lemma}{0}
\renewcommand{\thelemma}{A.\arabic{lemma}}

For simplicity, we denote $\|f\|_{Q,r}$ as the $L_r(Q)$-norm of
the function $f$. Let $T_n^\ast$ be a random vector composed of
$(X_1,\ldots,X_n)$ and $(W_{n1},\ldots,W_{nn})$. According to
\citep{h84}, we say that the conditional distribution of
$T_n^\ast$ given $\mathcal{X}_n$ converges weakly to some fixed
distribution $T$ in $P_X$-probability, denoted as
$``T_n^\ast\Longrightarrow T"$, if
\begin{eqnarray}
\sup_{f\in BL_1}\left|E_{W|\mathcal{X}_n}f(T_n^\ast)-\int
fdT\right|\overset{P_X}{\longrightarrow}0,\label{inter2}
\end{eqnarray}
where $BL_1$ is the class of Lispchitz functions bounded by $1$
and with Lipschitz norm $1$.

%

\subsection{Proof of Theorem~\ref{bconthm}}
Choose some $p''$ satisfying $p<p''<p'$. According to
Lemma~\ref{keylem} and the definition of $\rho_p$, it suffices to show that
$$\sup_{n}E_{XW}\|\sqrt{n}(\widehat\theta^\ast-\theta_0)\|^{p''}<\infty\;\;\;\mbox{and}\;\;\;\sup_{n}E_{X}\|\sqrt{n}(\widehat\theta-\theta_0)\|^{p''}<\infty.$$
The latter result is a special case of the former since we may
take $W_{ni}=1$ a.s. for $i=1,\ldots,n$. To show the former, it
suffices to show
\begin{eqnarray}
\sup_n
E_{XW}\{\sqrt{n}[\|\widehat\theta^\ast-\theta_0\|+d(\widehat\eta^\ast,\eta_0)]\}^{p''}<\infty.\label{finres}
\end{eqnarray}
To show (\ref{finres}), we need to partition the parameter space
$\Theta\times\mathcal{H}$ into ``shells'' $S_{j,n}$, i.e.,
\begin{eqnarray}
S_{j,n}=\left\{(\theta,\eta)\in\Theta\times\mathcal{H}:
2^{j-1}<\sqrt{n}(\|\theta-\theta_0\|+d(\eta,\eta_0))\leq
2^j\right\}\label{shell}
\end{eqnarray}
with $j$ ranging over integers, and then bound the
probability of each shell under Conditions M1-M2. For any fixed
$j_0>0$, we have
\begin{eqnarray*}
&&E_{XW}\{\sqrt{n}[\|\widehat\theta^\ast-\theta_0\|+d(\widehat\eta^\ast,\eta_0)]\}^{p''}\\
&\leq&
2^{(j_0-1)p''}P_{XW}(\sqrt{n}(\|\widehat\theta^\ast-\theta_0\|+d(\widehat\eta^\ast,\eta_0))\leq
2^{j_0-1})\\&&+\sum_{j=j_0}^\infty
2^{jp''}P_{XW}(2^{j-1}<\sqrt{n}(\|\widehat\theta^\ast-\theta_0\|+d(\widehat\eta^\ast,\eta_0))\leq 2^j)\\
&\leq&2^{(j_0-1)p''}+\sum_{j=j_0}^\infty
2^{jp''}P_{XW}\left(\sup_{(\theta,\eta)\in
S_{j,n}}\mathbb{P}_n^\ast\left( m(\theta, \eta)- m(\theta_0,
\eta_0)\right)\geq0\right)\\
&\leq&2^{(j_0-1)p''}+\sum_{j=j_0}^\infty
2^{jp''}P_{XW}\left(\sup_{(\theta,\eta)\in
S_{j,n}}(\mathbb{P}_n^\ast-P_X)\left( m(\theta, \eta)- m(\theta_0,
\eta_0)\right)\apgt\frac{2^{2j-2}}{n}\right),
\end{eqnarray*}
where the last inequality follows from Condition M1. By the
decomposition that
$(\mathbb{P}_n^\ast-P_X)f=n^{-1/2}(\mathbb{G}_n^\ast+\mathbb{G}_n)f$,
we can further bound the second term in the above by
\begin{eqnarray*}
&&\sum_{j=j_0}^\infty 2^{jp''}P_{XW}\left(\sup_{(\theta,\eta)\in
S_{j,n}}\mathbb{G}_n^\ast(m(\theta, \eta)-
m(\theta_0, \eta_0))\apgt\frac{2^{2j-2}}{\sqrt{n}}\right)\\
&+&\sum_{j=j_0}^\infty 2^{jp''}P_{X}\left(\sup_{(\theta,\eta)\in
S_{j,n}}\mathbb{G}_n(m(\theta, \eta)-
m(\theta_0, \eta_0))\apgt\frac{2^{2j-2}}{\sqrt{n}}\right)\\
&\leq&\sum_{j=j_0}^\infty
2^{jp''}\left[\left(\frac{2^j/\sqrt{n}}{2^{2j-2}/\sqrt{n}}\right)^{p'}+\left(\frac{2^j/\sqrt{n}}{2^{2j-2}/\sqrt{n}}\right)^{p'}
\right]\\
&\aplt&\sum_{j=j_0}^\infty 2^{j(p''-p')}
\end{eqnarray*}
The first inequality follows from Markov's inequality and
Condition M2. Now, we can conclude that
\begin{eqnarray*}
&&E_{XW}\{\sqrt{n}[\|\widehat\theta^\ast-\theta_0\|+d(\widehat\eta^\ast,\eta_0)]\}^{p'}\\
&\leq&2^{(j_0-1)p'}+\sum_{j=j_0}^\infty 2^{j(p''-p')}<\infty
\end{eqnarray*}
since we assume that $p''<p'$. This concludes the proof. $\Box$

\subsection{Conditional Slutsky's Lemma}\label{conslt}

Suppose $T_n^\ast\Longrightarrow T$
and $C_n\overset{P_X}{\longrightarrow} C$ for some vector $C$,
then we have
\begin{enumerate}
\item[(i)] $T_n^\ast+C_n\Longrightarrow T+C$;

\item[(ii)] $C_nT_n^\ast\Longrightarrow CT$;

\item[(iii)] $C_n^{-1}T_n^\ast\Longrightarrow C^{-1}T$ provided
$C\neq 0$,
\end{enumerate}
where $T_n^\ast$ and $C_n$ are random vectors composed of
$(X_1,\ldots,X_n)$ and $(W_{n1},\ldots,W_{nn})$, and
$(X_1,\ldots,X_n)$, respectively. In addition, the vector $C$ in (i) must be of the same dimension as $T$
and $C$ in (ii) \& (iii) can be a matrix.

{\it Proof:} Without loss of generality, we assume $C$ to be a
vector. If $C$ is a matrix, the conclusions in (ii) and (iii) are
still valid since the matrix multiplication and matrix inversion
are both continuous operations. We first show the conditional weak
convergence $(T_n^\ast, C_n)\Longrightarrow (T, C)$, and then
apply the conditional version of the continuous mapping Theorem,
i.e., Theorem 10.8 in \citep{k08}, to conclude the proof. We first
show the following result:
\begin{eqnarray}
\mbox{if}\;\;U_n^\ast\Longrightarrow
U\;\;\mbox{and}\;\;\|U_n^\ast-V_n^\ast\|\overset{P_X}{\longrightarrow}0,\;\;\mbox{then}\;V_n^\ast\Longrightarrow
U, \label{inter1}
\end{eqnarray}
where $U_n^\ast$ and $V_n^\ast$ are random vectors composed of
$(X_1,\ldots,X_n)$ and $(W_{n1},\ldots,W_{nn})$. For any $f\in
BL_1$, we have
\begin{eqnarray}
|E_{W|\mathcal{X}_n}f(U_n^\ast)-E_{W|\mathcal{X}_n}f(V_n^\ast)|\leq\epsilon
E_{W|\mathcal{X}_n}1\{\|U_n^\ast-V_n^\ast\|\leq\epsilon\}+2E_{W|\mathcal{X}_n}
1\{\|U_n^\ast-V_n^\ast\|>\epsilon\}\label{inter8}
\end{eqnarray}
for every $\epsilon>0$. The first term in the right hand side of
(\ref{inter8}) can be made arbitrarily small by choice of
$\epsilon$ while the second term converges to zero in
$P_X$-probability as $n\rightarrow\infty$. Thus, we claim
$$\sup_{f\in BL_1}\left|E_{W|\mathcal{X}_n}f(U_n^\ast)-E_{W|\mathcal{X}_n}f(V_n^\ast)\right|\overset{P_X}{\longrightarrow}0.$$
Considering the definition (\ref{inter2}) and
$U_n^\ast\Longrightarrow U$, we complete the proof of
(\ref{inter1}). According to (\ref{inter1}), it suffices to show
$(T_n^\ast, C)\Longrightarrow(T, C)$ since $\|(T_n^\ast,
C_n)-(T_n^\ast, C)\|=\|C_n-C\|\overset{P_X}{\longrightarrow}0$. It
is easy to show that for every bounded Lipschitz function
$(x,y)\mapsto f(x,y)$, the function $x\mapsto f(x,c)$ is also
bounded and Lipschitz continuous. Thus, if
$T_n^\ast\Longrightarrow T$, then we have
$$\sup_{f\in BL_1}\left|E_{W|\mathcal{X}_n}f(T_n^\ast, C)-Ef(T,
C)\right|\leq\sup_{f\in
BL_1}\left|E_{W|\mathcal{X}_n}f(T_n^\ast)-Ef(T)\right|\overset{P_X}{\longrightarrow}0.$$
Again, an application of (\ref{inter2}) completes the whole proof.
$\Box$

\subsection{An Inequality for $\|\cdot\|_{2,1}$-norm}\label{ineapp}
The following chain inequality is essentially Problem 2.9.1 of \cite{vw96}. We provide the proof for completeness.

For any $Y>0$ and $r>2$, we have
\begin{eqnarray}
\frac{1}{2}\|Y\|_2\leq\|Y\|_{2,1}\leq\frac{r}{r-2}\|Y\|_r,\label{inter7}
\end{eqnarray}
where $\|Y\|_r=(EY^r)^{1/r}$.

{\it Proof:} The first inequality is established as follows:
$$\|Y\|_2^2=2\int_0^\infty tP(|Y|>t)du=2\int_0^\infty t\sqrt{P(|Y|>t)}\sqrt{P(|Y|>t)}dt\leq2\|Y\|_{2,1}\|Y\|_2$$ by Markov's
inequality. For the second inequality, we have
\begin{eqnarray*}
\|Y\|_{2,1}&=&\left(\int_0^a+\int_a^\infty\right)\sqrt{P(|Y|>t)}dt\\&\leq&
a+\|Y\|_r^{r/2}\int_a^\infty
t^{-r/2}dt\\&\leq&a+\|Y\|_r^{r/2}\frac{2a^{1-r/2}}{r-2}\equiv U(a)
\end{eqnarray*}
for any $a>0$. It is easy to show that the minimal of $U(a)$ is
just $[r/(r-2)]\|Y\|_r$ when $a=\|Y\|_r$. This completes the proof
of the second inequality in (\ref{inter7}). $\Box$

\subsection{The $L_p$ Multiplier Inequality}\label{lpmin}
Let $W_{n}=(W_{n1},\ldots,W_{nn})'$ be non-negative exchangeable
random variables on $(\mathcal{W},\Omega,P_{W})$ such that, for
every $n$, $R_{n}=\int_{0}^{\infty}\sqrt{P_{W}(W_{n1}\geq
u)}du<\infty$. Let $Z_{ni}$, $i=1,2,\ldots,n$, be i.i.d. random
elements in
$(\mathcal{X}^{\infty},\mathcal{A}^{\infty},P_X^{\infty})$ with
values in $\ell^{\infty}(\mathcal{F}_{n})$, and write
$\|\cdot\|_{n}=\sup_{f\in\mathcal{F}_{n}}|Z_{ni}(f)|$. It is
assumed that $Z_{ni}$'s are independent of $W_{n}$. Then for any
$n_{0}$ such that $1\leq n_{0}<\infty$ and any $n>n_{0}$, the
following inequality holds for any $p\geq 1$:
\begin{eqnarray}
\left\|\left\|\frac{1}{\sqrt{n}}\sum_{i=1}^{n}
W_{ni}Z_{ni}\right\|_{n}\right\|_{P_{XW},p}&\leq&
n_{0}\left\|\left\|Z_{n1}\right\|_{n}\right\|_{P_X, p}\cdot\frac{
\left\|\max_{1\leq i\leq
n}W_{ni}\right\|_{P_W,p}}{\sqrt{n}}\nonumber\\&&+R_{n}^{1/p}\cdot\left\|\max_{n_{0}<
i\leq n}
\left\|\frac{1}{\sqrt{i}}\sum_{j=n_{0}+1}^{i}Z_{nj}\right\|_{n}\right\|_{P_X,
p}.\label{mulineq}
\end{eqnarray}

{\it Proof:} This Lemma generalizes the results in Lemma 4.1 of
\citep{wz96} where $p=1$. By the triangle inequality, we have
\begin{eqnarray*}
\left\|\left\|\frac{1}{\sqrt{n}}\sum_{i=1}^{n}
W_{ni}Z_{ni}\right\|_{n}\right\|_{P_{XW},p}&\leq&\left\|\left\|\frac{1}{\sqrt{n}}\sum_{i=1}^{n_0}
W_{ni}Z_{ni}\right\|_{n}\right\|_{P_{XW},p}+\left\|\left\|\frac{1}{\sqrt{n}}\sum_{i=n_0+1}^{n}
W_{ni}Z_{ni}\right\|_{n}\right\|_{P_{XW},p}.
\end{eqnarray*}
The first term in the above is trivially bounded by
$$n_{0}\left\|\left\|Z_{n1}\right\|_{n}\right\|_{P_X,
p}\cdot\frac{\left\|\max_{1\leq i\leq
n}W_{ni}\right\|_{P_W,p}}{\sqrt{n}}.$$ Denote $W_{n(i)}$ as the
$i$th ordered values of $W_{ni}$, i.e., $W_{n(1)}\geq
W_{n(2)}\geq\cdots\geq W_{n(n)}$. Note that $\|\|\sum_{i=1}^n
W_{ni}Z_{ni}\|_n\|_{P_{XW},p}=\|\|\sum_{i=1}^n
W_{n(i)}Z_{ni}\|_n\|_{P_{XW},p}$ since $W_n$ is assumed to be
exchangeable and $P_X^\infty$ is permutation invariant. We write
the second term as the following telescoping sum,
$$\sum_{i=n_0+1}^{n}W_{n(i)}Z_{ni}=\sum_{i=n_0+1}^n\sqrt{i}(W_{n(i)}-W_{n(i+1)})T_i,$$
where $T_i\equiv i^{-1/2}\sum_{j=n_0+1}^i Z_{nj}$ and
$W_{n(n+1)}\equiv 0$. Thus, we obtain that
\begin{eqnarray*}
\left\|\left\|\sum_{i=n_0+1}^{n}
W_{n(i)}Z_{ni}\right\|_{n}\right\|_{P_{XW},p}&\leq&\left\|\sum_{i=n_0+1}^n\sqrt{i}(W_{n(i)}-W_{n(i+1)})\|T_i\|_n\right\|_{P_{XW},p}\\
&\leq&\left\|\max_{n_0< i\leq
n}\|T_i\|_n\right\|_{P_X,p}\cdot\left\|\sum_{i=n_0+1}^n\sqrt{i}(W_{n(i)}-W_{n(i+1)})\right\|_{P_W,p}.
\end{eqnarray*}
Recalling the definition of $T_i$, it remains to show
\begin{eqnarray}
E_W\left(\sum_{i=n_0+1}^n\sqrt{i}(W_{n(i)}-W_{n(i+1)})\right)^{p}\leq
n^{p/2}R_n.\label{inter3}
\end{eqnarray}
Note that
\begin{eqnarray*}
&&\sum_{i=n_0+1}^n\sqrt{i}(W_{n(i)}-W_{n(i+1)})\\
&=&\sum_{i=n_0+1}^n\int_{W_{n(i+1)}}^{W_{n(i)}}\sqrt{i}du=\int_0^{W_{n(n_0+1)}}\sum_{i=n_0+1}^n
\sqrt{i}1\{W_{n(i)}\leq u\leq W_{n(i+1)}\}du\\&=&\int_0^{W_{n(n_0+1)}}\sqrt{\sharp\{r\geq 1: W_{n(r)}\geq u\}}du,
\end{eqnarray*}
which is bounded by
$$\left(\int_0^{W_{n(n_0+1)}}\left\{\sharp\{r\geq 1: W_{n(r)}\geq
u\}\right\}^{p/2}du\right)^{1/p}.$$
By taking the expectation $E_W$ in the above, we have shown that the left hand side of (\ref{inter3}) is bounded by
\begin{eqnarray*}
E_W\int_0^{W_{n(n_0+1)}}\left\{\sharp\{r\geq 1: W_{n(r)}\geq
u\}\right\}^{p/2}du&\leq& n^{p/2} E_W\int_0^{\infty}\left\{\sharp\{r\geq 1:
W_{n(r)}\geq u\}/n\right\}^{p/2}du\\
&\leq&n^{p/2} E_W\int_0^{\infty}\left\{\sharp\{r\geq 1:
W_{n(r)}\geq u\}/n\right\}^{1/2}du\\
&\leq&n^{p/2}\int_0^\infty\sqrt{P_W(W_{n1}\geq u)}du=n^{p/2}R_n
\end{eqnarray*}
based on the Jensen's inequality. This completes the whole proof.
$\Box$

\subsection{Verification of Condition (\ref{conmod4})}\label{conver}
Suppose that the $L_{p'}$ maximal inequality (\ref{conmod3}) and
bootstrap weight condition (\ref{boowei}) hold. If
$\|N_\delta\|_{P_X, p'}<\infty$, then we have Condition
(\ref{conmod4}) for each $p'\geq 1$.

{\it Proof:} We first apply the symmetrization argument to show
\begin{eqnarray}
\left\|\|\mathbb{G}_n^\ast\|_{\mathcal{N}_\delta}\right\|_{P_{XW},
p'}\aplt2\left\|\left\|\frac{1}{\sqrt{n}}\sum_{i=1}^{n}W_{ni}(\delta_{X_i}-P_X)\right\|_{\mathcal{N}_\delta}\right\|_{P_{XW},
p'}.\label{inter4}
\end{eqnarray}
Note that
\begin{eqnarray*}
\mathbb{G}_{n}^{\ast}=\frac{1}{\sqrt{n}}\sum_{i=1}^{n}(W_{ni}-1)\delta_{X_{i}}=\frac{1}{\sqrt{n}}\sum_{i=1}^{n}
(W_{ni}-1)(\delta_{X_{i}}-P_X)
\end{eqnarray*}
by Condition W2. Let $W_{n}'=(W_{n1}',\ldots,W_{nn}')$ be
exchangeable bootstrap weights generated from $P_{W'}$, an
independent copy of $P_{W}$. The bootstrap weight conditions $W1$
and $W2$ imply that $E_{W'}W_{ni}'=1$ for $i=1,\ldots,n$. Then, we
have
\begin{eqnarray*}
E_{XW}\|\mathbb{G}_n^\ast\|_{\mathcal{N}_\delta}^{p'}&=&
E_{XW}\left\|\frac{1}{\sqrt{n}}\sum_{i=1}^{n}(W_{ni}-1)(\delta_{X_i}-P_X)\right\|_{\mathcal{N}_\delta}^{p'}\\
&=&E_{XW}\left\|\frac{1}{\sqrt{n}}\sum_{i=1}^{n}(W_{ni}-E_{W'}W_{ni}')(\delta_{X_i}-P_X)\right\|_{\mathcal{N}_\delta}^{p'}\\
&\leq&E_{XW}E_{W'}\left\|\frac{1}{\sqrt{n}}\sum_{i=1}^{n}(W_{ni}-W_{ni}')(\delta_{X_i}-P_X)\right\|_{\mathcal{N}_\delta}^{p'}
\end{eqnarray*}
based on the Jensen's inequality and the reverse Fatou's Lemma. In
the end, a typical application of the symmetrization argument and
Minkowski's inequality concludes (\ref{inter4}).

To further bound the right hand side of (\ref{inter4}), we next
apply the $L_p$ multiplier inequality (\ref{mulineq}) with
$Z_{ni}=(\delta_{X_i}-P_X)$ and
$\mathcal{F}_n=\mathcal{N}_\delta$. This gives, due to Condition
W3,
$$\left\|\|\mathbb{G}_n^\ast\|_{\mathcal{N}_\delta}\right\|_{P_{XW},
p'}\aplt\left\|\left\|Z_{n1}\right\|_{\mathcal{N}_\delta}\right\|_{P_X,
p'}\cdot \frac{1}{\sqrt{n}}\left\|\max_{1\leq i\leq
n}W_{ni}\right\|_{P_W,p'}+\left\|\max_{n_{0}< i\leq n}
\left\|\frac{1}{\sqrt{i}}\sum_{j=n_{0}+1}^{i}Z_{nj}\right\|_{\mathcal{N}_\delta}\right\|_{P_X,
p'}$$ for any $1\leq n_{0}<\infty$ and $n>n_{0}$. For the last
term in the above, we can bound it by
\begin{eqnarray*}
&&\left\|\max_{n_{0}< i\leq
n}\left\|\frac{1}{\sqrt{i}}\sum_{j=1}^{i}
Z_{nj}\right\|_{\mathcal{N}_\delta}\right\|_{P_X,
p'}+\left\|\left\|\frac{1}{\sqrt{n_0}}\sum_{j=1}^{n_0}Z_{nj}\right\|_{\mathcal{N}_\delta}\right\|_{P_X,
p'}\\
&\leq&\left\|\max_{n_{0}< k\leq
n}\left\|\mathbb{G}_k\right\|_{\mathcal{N}_\delta}\right\|_{P_X,
p'}+\left\|\left\|\mathbb{G}_{n_0}\right\|_{\mathcal{N}_\delta}\right\|_{P_X,
p'}\\&\leq&2\left\|\max_{n_{0}\leq k\leq
n}\left\|\mathbb{G}_k\right\|_{\mathcal{N}_\delta}\right\|_{P_X,
p'}
\end{eqnarray*}
by the triangular inequality. In addition, we can bound
$\|\|Z_{n1}\|_{\mathcal{N}_\delta}\|_{P_X, p'}$ as
$$\|\|Z_{n1}\|_{\mathcal{N}_\delta}\|_{P_X, p'}=\|\|\delta_{X_1}-P_X\|_{\mathcal{N}_\delta}\|_{P_X, p'}\leq
\|\|\delta_{X_1}\|_{\mathcal{N}_\delta}\|_{P_X,
p'}+\|\|P_{X}\|_{\mathcal{N}_\delta}\|_{P_X, p'}\leq
2\|N_\delta\|_{P_X,p'}$$ due to the reverse Fatou's Lemma. Thus,
we obtain that
\begin{eqnarray*}
\left\|\|\mathbb{G}_n^\ast\|_{\mathcal{N}_\delta}\right\|_{P_{XW},
p'}&\aplt&
\|N_\delta\|_{P_X,p'}\cdot\frac{1}{\sqrt{n}}\left\|\max_{1\leq
i\leq n}W_{ni}\right\|_{P_W,p'}+\left\|\max_{n_{0}\leq k\leq
n}\left\|\mathbb{G}_k\right\|_{\mathcal{N}_\delta}\right\|_{P_X,
p'}\\
&\aplt&I+II.\nonumber
\end{eqnarray*}

Considering Condition (\ref{boowei}) and Lemma 4.7 of
\citep{pw93}, we have $n^{-1/2}E_W(\max_{1\leq i\leq
n}W_{ni}^{p'})\longrightarrow 0$. The inequality that
$\|\max_{1\leq i\leq n}W_{ni}\|_{P_W, p'}\leq E_W(\max_{1\leq
i\leq n}W_{ni}^{p'})$ (due to $\max_{1\leq i\leq n}W_{ni}^{p'}\geq
1$) implies
\begin{eqnarray}
\frac{1}{\sqrt{n}}\left\|\max_{1\leq i\leq
n}W_{ni}\right\|_{P_W,p'}=o(1).\label{wo1}
\end{eqnarray}
Since $\|N_\delta\|_{P_X, p'}$ is assumed to be finite, the above
term $I$ converges to zero, and thus is smaller than arbitrary
$\delta>0$ for sufficiently large $n$. For any positive r.v. $Y$,
it is easy to prove that
$$EY^q=\int_0^{\infty}qt^{q-1}P(Y>t)dt\;\;\;\mbox{for any}\;q>0.$$
The L\'{e}vy's inequality, i.e., Proposition A.1.2 in
\citep{vw96}, implies that
$$P\left(\max_{k\leq n}\|\mathbb{G}_k\|_{\mathcal{N}_\delta}>\lambda\right)\leq 2P\left(\|\mathbb{G}_n\|_{\mathcal{N}_\delta}>\lambda\right)\;\;\;\mbox{for every}\;\lambda>0.$$
Thus, we have that
$II\leq2^{1/p'}\|\|\mathbb{G}_n\|_{\mathcal{N}_\delta}\|_{P_X,
p'}$. This concludes the whole proof. $\Box$

\end{document}